\documentstyle[]{article}
\def\a{\alpha}
\def\b{\beta}
\def\g{{\bf g}}
\def\bb{{\bf b}}
\def\h{{\bf h}}
\def\G{\Gamma}
\def\D{\Delta}

\def\l{\lambda}
\def\L{\Lambda}
\def\m{\mu}

\def\Aa{{\cal A}}
\def\A*{{\cal A}^*}

\def\ea{e_\a}
\def\e-a{e_{-\a}}

\newcommand{\bn}{\begin{equation}}
\newcommand{\ed}{\end{equation}}

\newtheorem{proposition}{Proposition}[section]

\newtheorem{corollary}{Corollary}[section]
\newcommand{\hsp}{\mbox{$\hspace{.3in}$}}

\begin{document}
%%%%%Subject: hirota1
%\begin{titlepage}
\begin{center}
{\LARGE\bf
 Intertwining operators and \\ Hirota bilinear equations}\\
\vspace{.5cm}
{talk given at the International Conference ''Modern problems of \\
 Quantum Field  Theory, Quantum Gravity and Strings'',\\
 Alushta, June 10-20, 1994}\\
\vspace{1cm}
{\large\bf S. Kharchev\footnote{E-mail: kharchev@lpi.ac.ru}}
\vspace{.5cm}\\
  Institute of
Theoretical and Experimental Physics \\ 117259 Moscow, Russia,
\vspace{.2cm}\\
{\large and}
\vspace{.2cm}\\
Theory Department ,  P.N.Lebedev Physics
Institute,\\ Leninsky prospect 53, Moscow,~117~924, Russia
\vspace{.5cm}\\
{\large\bf  S. Khoroshkin\footnote{E-mail:
 khor@s43.msk.su, khor@mpim-bonn.mpg.de},
%%\footnote{Partially supported by
%% ISF grant MBI 000},
D.Lebedev\footnote{E-mail: lebedev@vxdesy.desy.de}}
\vspace{.5cm}
 \\ Institute of
Theoretical and Experimental Physics \\ 117259 Moscow, Russia,
\end{center}
\date{}
\vspace{1.5cm}

%%\maketitle
\begin{abstract}
We give an interpretation of Hirota relations for $\tau$-functions of
  hierarchies of integrable equations in terms of intertwining operators.
This gives possibility to generalize the relations to the case of finite-
dimensional Lie algebras and quantized universal enveloping algebras.
 An example of $U_q(sl_2)$ is presented.
\end{abstract}
%\end{titlepage}
\setcounter{section}{0}
\section{Motivations. $\tau$-function and matrix elements.}
One of the most popular languages for algebraic study of nonlinear integrable
 equations is the language of $\tau$-function. Usually $\tau$-function is
 connected with initial variable $u$ via some logarithmic derivatives and can
 be treated as a point of infinite-dimensional grassmanian or of other
 homogeneous space for certain infinite-dimensional group $G$. For instance,
 for KP hierarchy  $u=  \frac{\partial^2\log \tau}{\partial x_1^2}$ and
$\tau$ is a point of $GL_\infty$ -orbit $X$ of highest weight vector
 $v_0$ in basic representation $\Lambda_0$ of $GL_\infty$ in
 semiinfinite forms. The coordinates on $X$ are given by boson-fermion
 correspondence and the equations of motion in Hirota bilinear form
 could be obtained by Kac--Wakimoto construction \cite{Kac}, \cite{KW}
from quadratic Plucker-type equations of the orbit
 $X$. These equations are
$$ \Omega \tau\otimes\tau =c\cdot \tau\otimes\tau$$
where $\Omega$ is an operator commuting with $g\otimes g$, $g\in G$ , $c$ is a
constant which one can compute acting by $\Omega$
 on tensor product of highest weight vectors:
$$ \Omega v_0\otimes v_0=c\cdot v_0\otimes v_0$$
For KP hierarchy $\Omega= \sum_{j\in{\bf Z}}\psi_j\otimes\psi_j^*$, where
 $\psi_j$ and $\psi_j^*$ are free fermions. For hierarchies arising from
 level one irreducible representations of a Kac--Moody algebra $\widehat{\g}$
 an element $\Omega$ is divided Casimir operator, representing
Killing form for
 $\widehat{\g}$:
$$\Omega =\sum e_i\otimes e^i$$.

Equivalently, $\tau$-function for KP hierarchy can be defined as a matrix
 element \cite{D}
\bn
\tau^g(x_1,x_2,\ldots)=<v_0\mid e^{H(\overline{x})}g\mid v_0>
\label{M1}
\ed
where $g\in GL_\infty$ and
$$e^{H(\overline{x})} = \exp(\sum_{i=1}^\infty x_ka_k)$$
is a flow defined by commutative subalgebra $\{ a_k, k\geq 1\}$ of Lie algebra
 $gl_\infty$:
$$ a_k=\sum_{j\in{\bf Z}}\psi_j\psi_{j+k}^*$$
or by a commutative part of appropriate Heisenberg subalgebra for a Kac--Moody
 algebra $\widehat{\g}$.

In these notations the basic relation which produces Hirota equations after
 the expansion over the diagonal is again the main property of $\Omega$:
 \bn
 \Omega g\otimes g =g\otimes g \Omega
\label{M2}
\ed
 For instance, in KP case we have due to (\ref{M2}):
$$\sum_{j\in{\bf Z}}<v_1\mid e^{H(\overline{x})}\psi_jg\mid v_0>
<v_{-1}\mid e^{H(\overline{y})}\psi_j^*g\mid v_0>=$$
\bn
=\sum_{j\in{\bf Z}}<v_1\mid e^{H(\overline{x})}g\psi_j\mid v_0>
<v_{-1}\mid e^{H(\overline{y})}g\psi_j^*\mid v_0>=0
\label{M3}
\ed
 where $<v_k\mid$ is a left vacuum of charge $k$:
$$<v_k\mid = <v_0\mid \psi_0^*\cdots \psi_{k-1}^*,\hsp k\geq 0,$$
$$<v_k\mid = <v_0\mid \psi_{-1}\cdots \psi_{-k},\hsp k< 0,$$
 and, moving fermions $\psi_j$ and $\psi_j^*$ to the left in LHS of
(\ref{M3}) we get desired functional relation in Hirota form
\bn
\sum_{j\in{\bf Z}}S_j(2\overline{y})S_{j+1}
(-\widetilde{\partial}_{{y}})
\tau^g(\overline{x}+\overline{y})\tau^g(\overline{x}-\overline{y})=0
\label{M4}
\ed
where $\widetilde{\partial}_{{y}} = \partial_{y_1}, \frac{1}{2}
\partial_{y_2}, \ldots$ and $S_j$ are Schur polynomials.

%%%%%%%%%Subject: hirota2

\setcounter{equation}{0}
\section{Hirota relations and Intertwining operators}
We can see that the procedure of getting the bilinear equations for
$\tau$-function could be broken into two steps:

{\it (i)} To find some (commutative) algebra which bosonize the representation;

{\it (ii)} To find bilinear relations for matrix elements of representation.

We could also have $\tau$-functions for different representations of the same
$G$. If we want to have corresponding $\tau$-functions being connected by a
 system of equations we should impose that all the representations involved
could be bosonized  by the same commutative algebra.

There is a general simple procedure of getting bilinear relations on matrix
 elements. Let us denote the matrix element $<ugv>$ of $g$ in representation
 $V$ by $\tau^g_{u,v}$:
\bn
\tau^g_{u,v}=<ugv>,\hsp v\in V,\;\; u\in V^*,\;\; g\in G
\label{H1}
\ed

Let now $V_1, V_2, V_3, V_4$ be representations of $G$, $v_i\in V_i$, $u_i\in
V_i^*$ be some fixed vectors and $\G = \sum_i\G^1_i \otimes \G^2_i$ be an
intertwining operator
\bn
\G :V_1\otimes V_2 \rightarrow V_3\otimes V_4
\label{H2}
\ed
that is,
\bn
\G \D(x)=\D(x)\G
\label{H3}
\ed
where $\D$ is comultiplication in $G$.
Then we have, since $\D(g)=g\otimes g$ for $g\in G$ and by (\ref{H3}),
 the following
\begin{proposition}
\bn
\sum_i\tau^g_{u_3\G^1_i,v_1}\tau^g_{u_4\G^2_i,v_2}=
\sum_i\tau^g_{u_3,\G^1_iv_1}\tau^g_{u_4,\G^2_iv_2}
\label{H4}
\ed
\label{prop1}
\end{proposition}
%%In other words, the symbol
%%$$\tau^g_{u_3\G^1v_1}\tau^g_{u_4\G^2v_2}$$
%%is correctly defined.
We should also keep bosonization picture for a tensor $\G$ if it is constructed
 from some ''local'' data. There is a description of such a construction.

Let us fix some ''vector'' representation $W$ of Lie algebra $\g$. It is the
 first fundamental representation of $\g$ for finite-dimensional simple $\g$;
 vector representation for $gl_\infty$ and  evaluation
 representation for affine algebra $\widehat{\g}$, constructed by affinization
 of the first fundamental representation of $\g$.

Let $$\Phi_W: V_\l\otimes W\rightarrow V_\m\hsp{\mbox{and}}\hsp
\Psi^W:V_{\l'}\rightarrow W\otimes V_{\m'}$$
be intertwining operators (that is $\Psi^Wx = \D(x)\Psi^W$,
 $x\Phi_W = \Phi_W\D(x)$). We call them, as in \cite{XXZ}, vertex operators.

Then their composition
\bn
\G=(\Phi_W\otimes id)(id\otimes\Psi^W):V_\l\otimes V_{\l'}\rightarrow
V_\m\otimes V_{\m'}
\label{H5}
\ed
is an intertwining operator. Here $V_\l, V_{\l'}, V_\m, V_{\m'}$ are
representations of $\g$ with  highest weights $\l, \l', \mu, \m'$.

We can introduce the components of intertwining operators $\Phi_{W,i}$
and $\Psi^{W,i}$ , where $w_i$ is a basis of $W$ as follows:  $$
\Phi_{W,i}(v)=\Phi_W(v\otimes w_i),$$ $$ \Psi^W(v)=\sum_iw_i\otimes
\Psi^{W,i}(v).$$ Then $\G=\sum_i \Phi_{W,i}\otimes \Psi^{W,i}$.
\vspace{.3cm}

{\bf Example 1}. Let $W$ be a vector representation of $gl_\infty$, $e_{i,j}
(w_j) = w_i$ ($W$ could be realized in a space of free fermions $\psi_j$),
 $\Lambda_k$ be a representation of $gl_\infty$ in semiinfinite forms of
 charge $k$. Let
 $$\Phi_W: \Lambda_k\otimes W\rightarrow \Lambda_{k+1}, \hsp
   \Psi^W: \Lambda_k\rightarrow W\otimes\Lambda_{k-1}$$
be corresponding vertex operators. Their components $\Phi_{W,i}$ and
$\Psi^{W,i}$ are known as operations of insertion and of removing a
 fermion $\psi_i$. If we apply tensor $\G$ , constructed by the relation
 (\ref{H5}) to a tensor product $\L_0\otimes\L_0$:
$$\G:\L_0\otimes\L_0\rightarrow \L_1\otimes\L_{-1}$$
then the equations (\ref{H4}) turn to be Hirota relations (\ref{M4}) for KP
 hierarchy if we put, according to (\ref{M3}), $u_3= <v_1
 \mid e^{H(\overline{x})}$,
 $u_4 = <v_{-1}\mid e^{H(\overline{y})}$, $v_1= v_2 = \mid v_0>$.
  Moreover, following \cite{UT}, \cite{D}, we can input an evolution of $\tau$
 by ''negative'' times $y_j$ and define $\tau_n^g(\overline{x},\overline{y})$
 as
$$\tau_n^g(\overline{x},\overline{u})=
<v_n\mid e^{H(\overline{x})}g e^{H'(\overline{u})}\mid v_n>$$
where $H'(\overline{u}) = \sum_{k\geq 1}u_ka_{-k}$, $a_{-k}=
\sum_{j\in{\bf Z}}\psi_j\psi_{j-k}^*$ and $v_n$ is a vacuum of charge $n$.

Then the equation (\ref{M4}) for
$$\G:\L_n\otimes\L_m\rightarrow \L_{n+1}\otimes\L_{m-1}$$
turns to Hirota relations
$$
\sum_{j\in{\bf Z}}S_j(2\overline{y})S_{j+n-m+1}
(-\widetilde{\partial}_{{y}})
\tau^g_n(\overline{x}+\overline{y},\overline{u}+\overline{v})
\tau^g_m(\overline{x}-\overline{y},\overline{u}-\overline{v})=$$
\bn
\sum_{j\in{\bf Z}}S_{j+n-m+1}(-2\overline{v})S_{j}
(\widetilde{\partial}_{{v}})
\tau^g_{n+1}(\overline{x}+\overline{y},\overline{u}+\overline{v})
\tau^g_{m-1}(\overline{x}-\overline{y},\overline{u}-\overline{v})
\label{H6}
\ed

Technically for deriving (\ref{H6}) one needs only the commutation relations
 of vertex operators $\Phi_W$ and $\Psi^W$ with evolution operators
 $e^{H(\overline{x})}$ and $ e^{H'(\overline{u})}$ and the action of
 $\Phi_{W,i}$ and $\Psi^{W,i}$ on highest weight vectors. If we put
$$ \Psi(z)= \sum_{i\in{\bf Z}}\Phi_{W,i}z^i, \hsp
   \Psi^*(z)= \sum_{i\in{\bf Z}}\Psi^{W,i}z^{-i}$$
then these relations are standard ones \cite{D}:
$$e^{H(\overline{x})}\Psi(z)= e^{\sum_{n\geq 1}x_nz^n}\Psi(z)
e^{H(\overline{x})}, \hsp
e^{H'(\overline{u})}\Psi(z)= e^{\sum_{n\geq 1}u_nz^{-n}}\Psi(z)
e^{H'(\overline{u})},$$
$$\Psi^*(z)e^{H(\overline{x})}=e^{\sum_{n\geq 1}x_nz^n}
e^{H(\overline{x})}\Psi^*(z), \hsp
\Psi^*(z)e^{H'(\overline{u})}=e^{\sum_{n\geq 1}u_nz^{-n}}
e^{H'(\overline{u})}\Psi^*(z),$$
$$\Psi(z)\mid v_n>= z^{n+1}e^{j'(z)}\mid v_{n+1}>,\hsp
  <v_n\mid \Psi(z)=z^n<v_{n-1}\mid e^{-j(z)},$$
$$\Psi^*(z)\mid v_n>= z^{-n}e^{-j'(z)}\mid v_{n-1}>,\hsp
  <v_n\mid \Psi^*(z)=z^{-(n+1)}<v_{n+1}\mid e^{j(z)}$$
where
$$j(z) =\sum_{n\geq 1} a_n\frac{z^{-n}}{n},\hsp
 j'(z) =\sum_{n\geq 1} a_{-n}\frac{z^{n}}{n}$$.

For KdV and NLS hierarchies one should use level one irreducible
 representations $\L_0$ and $\L_1$ of $\widehat{sl}_2$  bosonized
 via correspondingly principal and homogenious Heisenberg subalgebras, and
 with ''vector'' representation $W_z$ being two-dimensional evaluation
 representation (in principal or homogeneous gradation),
\bn
 \G = \oint \frac{dz}{z}\Phi_{W,a}(z)\otimes\Psi^{W,a}(z).
\label{integral}
\ed
where $a =\pm$ is an index of two-dimensional representation $W_z$ of
$\widehat{sl}_2$.

We should like to point out that in the case of affine algebra $\widehat{\g}$
vertex operators $\Phi_W$ and $\Psi^W$ are actually intertwining operators
 for algebra $\widehat{\g}'$ where the full algebra $\widehat{\g}$ is obtained
 by adding grading element $d$ to $\widehat{\g}'$. The integral in
 (\ref{integral}) arise since $\G$ should be intertwining operator for
 the full algebra $\widehat{\g}$.

{\bf Example 2}. Using the technique of vertex operators we can describe
Hirota presentation of finite-dimensional Toda lattice.

Let, for instance, ${\g} = A_n;\; \omega_1, \ldots ,\omega_n$ be
fundamental representations of $A_n$, and
$$\tau_n^g(\overline{x},\overline{u})=
<v_n\mid e^{H(\overline{x})}g e^{H'(\overline{u})}\mid v_n>$$
where $v_n$ is a highest weight vector of $\omega_n$,
$$H(x)=\sum_{k=1}^nx_kI_k, \hsp H'(u)=\sum_{k=1}^nu_kI_{-k},$$
$$ I_k=\sum_{i=1}^{n+1-k}e_{i,i+k},\hsp I_{-k}=\sum_{i=1}^{n+1-k}
e_{i+k,i},$$
 $e_{i,k}\in sl(n+1)$ are matrix units, $W = \omega_1$.
If we use vertex operators
$$\Phi_{\omega_1}: \omega_k\otimes \omega_1\rightarrow \omega_{k+1}, \hsp
   \Psi^{\omega_1}: \omega_k\rightarrow \omega_1\otimes\omega_{k-1}$$
and, just as in example 1, introduce $\G: \omega_k \otimes\omega_m
\rightarrow \omega_{k+1} \otimes\omega_{m-1}$
then we get the relations for $\tau_n, \tau_m, \tau_{n+1}$ and $\tau_{m-1}$
 analogous to (\ref{H6}) (actually, higher times in (\ref{H6}) should be
 freezed in appropriate way).

%%%%%Subject: hirota3

\setcounter{equation}{0}
\section{Noncommutative variant}
The basic relation (\ref{H4}) on matrix elements of  representations could be
 rewritten as well for matrix elements of quantum deformation $U_q(\g)$ of
 enveloping algebra $U(\g)$, or, more generally, for matrix elements of
 representations of a Hopf algebra $\Aa$. For this we should treat matrix
 elements as elements of dual Hopf algebra $\A*$ ($\A* = Fun_q(G)$ if $\Aa =
 U_q(\g)$).

Let us introduce some notations.

Let $\Aa$ be a Hopf algebra, $\A*$ be a dual Hopf algebra and $\{ ,\}: \A*
\otimes \Aa \rightarrow {\bf C}$ be canonical Hopf pairing. Let $V$ be a
 representation of algebra $\Aa$, $v\in V, u\in V^*$. Define an element
 $\tau_{u,v}$ by the relation
\bn
\{ \tau_{u,v},x\}= <u\mid x \mid v>,\hsp x\in \Aa
\label{Q1}
\ed
Let now $V_1, V_2, V_3, V_4$ be representations of $\Aa$, $v_i\in V_i$, $u_i\in
V_i^*$ be some fixed vectors and $\G = \sum_i\G^1_i \otimes \G^2_i$ be an
intertwining operator
\bn
\G :V_1\otimes V_2 \rightarrow V_3\otimes V_4
\label{Q2}
\ed
that is,
\bn
\G \D(x)=\D(x)\G
\label{Q3}
\ed
\begin{proposition}
\bn
\sum_i\tau_{u_3\G^1_i,v_1}\tau_{u_4\G^2_i,v_2}=
\sum_i\tau_{u_3,\G^1_iv_1}\tau_{u_4,\G^2_iv_2}
\label{Q4}
\ed
\label{th1}
\end{proposition}
%%In other words, the symbol
%%$$\tau_{u_3\G^1v_1}\tau_{u_4\G^2v_2}$$
%%is correctly defined.

{\bf Proof}. We should prove that $\{ {\mbox{LHS of (\ref{H4})}},x\} =
\{ {\mbox{RHS of (\ref{H4})}},x\}$ for any $x \in \Aa$.
Let us remind that for a
Hopf pairing $\{ ,\}$ we have, in particular $\{ ab,c\} =\{ a\otimes b,\D(c)
\}$ , or using the notation $\D(c)=c^{(1)}\otimes c^{(2)}$,
$\{ ab,c\} =\{ a,c^{(1)}\}\{ b,c^{(2)}\}$.

 We have
$$\sum_i \{\tau_{u_3\G^1_i,v_1}\tau_{u_4\G^2_i,v_2},x\}=
\sum_i \{\tau_{u_3\G^1_i,v_1}\otimes\tau_{u_4\G^2_i,v_2},\D(x)\}=$$
$$= \sum_i <u_3\G^1_i\mid x^{(1)} \mid v_1><u_4\G^2_i\mid x^{(2)} \mid v_2>=
<u_3\otimes u_4 \G \mid\D(x) \mid v_1\otimes v_2>=$$
$$=<u_3\otimes u_4 \mid \D(x)\G \mid v_1\otimes v_2>=
\sum_i <u_3\mid x^{(1)}\mid\G^1_iv_1><u_4 \mid x^{(2)}\mid\G^2_i v_2>=$$
$$=\sum_i \{\tau_{u_3,\G^1_iv_1}\otimes\tau_{u_4,\G^2_iv_2},\D(x)\} =
\sum_i \{\tau_{u_3,\G^1_iv_1}\tau_{u_4,\G^2_iv_2},x\}.$$
\vspace{.2cm}

Let now $g: \A* \rightarrow $End$\, Y$ be a representation of an algebra
$\A*$.  Let us denote by $\tau_{u,v}^g$ an operator $ g\tau_{u,v} \in
 $End$\, Y$. We have immediately the following generalization of
Proposition \ref{prop1}.
\begin{corollary} In the setting of Theorem
\ref{th1} we have an equality \bn
\sum_i\tau^g_{u_3\G^1_i,v_1}\tau^g_{u_4\G^2_i,v_2}=
\sum_i\tau^g_{u_3,\G^1_iv_1}\tau^g_{u_4,\G^2_iv_2}
\label{Q5}
\ed
\label{corollary}
\end{corollary}
In classical limit $q=1$ the Hopf pairing $\{ ,\} : Fun(G)
\otimes U(\g)\rightarrow {\bf C}$ is given by the relation
\bn
\{ f,x\} =L_xf(e)\hsp {\mbox{for}}\;\; f(a)\in Fun(G),\;\;\; a\in G.
\label{fx}
\ed
where $L_x$ is rightinvariant differential operator on $G$
corresponding to element $x \in U(\g)$.
 Comparing (\ref{fx}) with (\ref{Q1}) we see that
$\tau_{u,v}$ can be identified with the following function of $a \in G$:
$ \tau_{u,v}= <u\mid a\mid v>.$
A representation $g$ of algebra $Fun(G)$ is an evaluation of a function in
 a point $g\in G$, that is
$$ \tau_{u,v}^g=\tau_{u,v}(g)=<u\mid g\mid v>.$$
Therefore, the classical limit of (\ref{Q5}) coincides with (\ref{H4}).

We have the same construction of  tensor $\G$ via vertex operators in
 noncommutative  case:\\ for
$$\Phi_W: V_\l\otimes W\rightarrow V_\m\hsp{\mbox{and}}\hsp
\Psi^W:V_{\l'}\rightarrow W\otimes V_{\m'}$$
being intertwining operators (that is $x\Phi_W = \Phi_W\D(x)$ and
 $\Psi^Wx = \D(x)\Psi^W$ we have
$$\G=(\Phi_W\otimes id)(id\otimes\Psi^W):V_\l\otimes V_{\l'}\rightarrow
V_\m\otimes V_{\m'}.$$
Here $V_\l, V_{\l'}, V_\m, V_{\m'}$ are representations of $\Aa$; for
 $\Aa = U_q(\g)$ they are highest weight representations.
Note that it is important now that $\Psi^W$ creates $W$ to the left and
 $\Phi_W$ annihilate $W$ from the right, that is, $\Phi_W$ is type I vertex
 operator and $\Psi^W$ is type II vertex operator in terminology of
 \cite{XXZ}, \cite{JM}.

Alternatively, we can use dual vertex operators
$$\Phi^W:V_\l\rightarrow V_\m\otimes W\hsp{\mbox{and}}\hsp
\Psi_W:W\otimes V_{\l'}\rightarrow V_{\m'}$$
Then
$$\G=(id \otimes \Psi_W)(\Phi^W\otimes id):V_\l\otimes V_{\l'}\rightarrow
V_\m\otimes V_{\m'}$$
is again an intertwining operator constructed from type I and type II
 intertwining operators. The two presentations of vertex operators and of
tensor
 $\G$ are actually equivalent to replacement of $W$ by $W^{*S}$ where $S$ is an
 antipode in $\Aa$ due to canonical isomorphisms of the components of
vertex operators \cite{JM}:
 $$\Phi^{W,i}=
\Phi_{W^{*S'},i}:V_\l\rightarrow V_\m,$$ $$\Psi^{W,i}=
\Psi_{W^{*S},i}:V_{\l'}\rightarrow V_{\m'}$$ Here the dual bases $\{
w_i\}$ and $\{ w^i\}$ of $W$ and $W^*$ are used; $S'$ is an inverse to
 $S$.

The defining commutation relations
\bn
\Phi^Wx = \D(x)\Phi^W,\hsp x\Phi_W = \Phi_W\D(x),
\label{Q6}
\ed
\bn
\Psi^Wx = \D(x)\Psi^W,\hsp x\Psi_W = \Psi_W\D(x),
\label{Q7}
\ed
could be  rewritten also in terms of the components of the vertex operators
 (compare with \cite{FR}). We summarize them in the following proposition.
\begin{proposition}
The components of vertex operators satisfy the relations
$$ S(x^{(1)})\Phi^{W,i}x^{(2)}=\sum_j\rho_j^i(x)\Phi^{W,j},\hsp
 S'(x^{(2)})\Psi^{W,i}x^{(1)}=\sum_j\rho_j^i(x)\Psi^{W,j},$$
 $$x^{(2)}\Phi_{W,i}S'(x^{(1)})=\sum_j\rho_i^j(x)\Phi_{W,j},\hsp
x^{(1)}\Psi_{W,i}S(x^{(2)})=\sum_j\rho_i^j(x)\Psi_{W,j}.$$
Here $\rho_i^j$ is a matrix for the action of $x\in\Aa$  in $W$:
$$xw_i=\sum_j\rho_i^j(x)w_j.$$
\label{Kharchev}
\end{proposition}
If the Hopf algebra $\Aa$ is quasitriangular, then the components of the vertex
 operators satisfy also the commutation relations
 of Zamolodchikov--Faddev  algebra \cite{ZZ}, \cite{F}, \cite{JM}.
We do not use them here.

Let now $\Aa$ be a quantized universal algebra $U_q(\g)$. In this case an
 $R$-matrix gives one more example of  tensor $\G$. Let us look to
 the basic equation (\ref{Q4}) in this particular case. Let $v_i$ be a basis
 of a representation $V$ of $U_q(\g)$, $u_i$ be a basis
 of a representation $U$, $<v_i,v^j> =<u_i,u^j> =\delta_{i,j}$, and
 $\widehat{R}:V\otimes U\rightarrow U\otimes V$ be an intertwining operator,
$$\widehat{R}(v_i\otimes u_j)= \sum_{a,b}R_{i,j}^{a,b}u_a\otimes v_b$$.
Denote $T_i^j =\tau_{v^j,v_i}$, ${T'}_i^j =\tau_{u^j,u_i}$, $T= \{
 T_i^j\}$, $T'= \{ {T'}_i^j\}$. Then the equation (\ref{Q4}) means that
$$ R_{i,j}^{a,b}T_a^k{T'}_b^l=R_{c,d}^{k,l}{T'}^c_iT_j^d$$
which is a traditional equation \cite{FRT} for matrix elements
\bn
RT_1{T'}_2={T'}_2T_1R,\hsp R=P\widehat{R}.
\label{RTT}
\ed
 Unfortunately, an $R$-matrix does not fit for deducing quantum Hirota-type
 equations; it is not built from local blocks and we cannot present in general
 a bosonized form for $R$.

We would like to mention that the equations (\ref{Q4}) with $\G$ being
constructed from
 vertex operators for a fixed vector representation $W$ give an alternative to
 (\ref{RTT}) description of an algebra $Fun_q(G)$ of matrix elements: there
 are no $RTT$-relations and we never permute  matrices $T$ and $T'$, but,
 unlike to (\ref{RTT}), where two representations of $U_q(\g)$  $V$ and $U$
 are involved, we have four representations $V_i$ which participate an equation
(\ref{Q4}). One may ask whether these equations describe completely an
algebra $Fun_q(G)$. Moreover, there are two variants of the question.

First, we can consider matrix elements of all integrable representations of
$U_q(\g)$ with all the equations (\ref{Q4}) built from vertex operators
 $\Phi^W$, $\Phi_W$, $\Psi^W$ and $\Psi_W$ for a given $W$. On the other hand,
 we may restrict ourselves by matrix elements of fundamental representations
 and  by vertex operators whose components act between them. In the case  of
 ${\g} =A_n$ one can show that both algebras are isomorphic to $Fun_q(G)$. It
will be
 interesting to find a proof of analogous statement in general case.
%%%%%%%%%%%%%%%%%%%%%%%%%%%%%%%
%%%%%%%%%%%%%%%%%%%%%%%%%%%%%%%
\setcounter{equation}{0}
\section{An example of $U_q(sl_2)$: $q$-diference Liouville equation}
Here we  present a toy example of $U_q(sl_2)$ \cite{GKLMM}.
An evolution in a representation of arbitrary spin $j$ in this case can be
 defined via $q$-exponents of generators of $U_q(sl_2)$. As a consequence of
 basic equations (\ref{Q4}) we get noncommutative $q$-difference analog
 of Liouville
 equation and its generalization to arbitrary spin $j$. Another approach
  is deveped in \cite{GS}.

Let us fix the notations. We use the following description of $U_q(sl_2)$:
$$k_\a e_{\pm\a}=q^{\pm 2}e_{\pm\a}k_\a, \hsp
 [\ea ,\e-a ]=\frac{k_\a -k_\a^{-1}}{q-q^{-1}},\hsp k_\a =q^{h_\a},$$
 $$ \D(\ea)=\ea \otimes 1+k_\a^{-1}\otimes\ea ,\hsp
 \D(\e-a)=1\otimes\e-a +\e-a\otimes k_\a ,$$
 $$\exp_q(x)= 1+x+\ldots +\frac{x^n}{(n)_q!}+\ldots$$
where
$$ (n)_q=\frac{q^n-1}{q-1}\hsp {\mbox{and}}\hsp [n]_q=
 \frac{q^n-q^{-n}}{q-q^{-1}}.$$

For a representation $V_j$ of spin $j$ we define $\tau_j(t,s) \in
Fun_q(G)$ as
 $$\tau_j(t,s)=\tau_{<j\mid\exp_{q^2}(t\ea ),
\exp_{q^{-2}}(s\e-a)\mid j>}.$$
Equivalently, the Hopf pairing of $\tau_j(t,s)$ with $x \in U_q(\g)$  is
 given by the relation
$$\{ \tau_j(t,s), x\}= <j\mid\exp_{q^2}(t\ea )\mid
\; x\;\mid\exp_{q^{-2}}(s\e-a)\mid j>$$
where $\mid j> $ is a
highest weight vector, $k_\a \mid j> = q^{2j}\mid j>$.  We denote by
$w_+$ and $w_-$ the standard basis of two-dimensional representation $W$
 of $U_q(sl_2)$: $\ea w_-=w_+$, $\e-a w_+=w_-$ and by $\Phi_\pm$  and
$\Psi^\pm$ the corresponding components of vertex operators $\Phi_W$ and
 $\Psi^W$:  $$\Phi_\pm: V_j\rightarrow V_{j'},\hsp \Psi^\pm:
 V_j\rightarrow V_{j''}$$ We can substitute $\exp_{q^2}(t\ea )$ or
$\exp_{q^{-2}}(s\e-a)$ instead of $x$ into defining relations $$
x\Phi_W=\Phi_W\D(x), \hsp \Psi^Wx=\D(x)\Psi^W.$$ The elements
$\exp_{q^2}(t\ea )$ and $\exp_{q^{-2}}(s\e-a)$ satisfy the following
factorization properties  with respect to comultiplication
\bn
\D(\exp_{q^2}(t\ea ))=\exp_{q^2}(tk_\a^{-1}\otimes\ea )
\exp_{q^2}(t\ea \otimes 1),
\label{sl1}
\ed
 \bn
\D(exp_{q^{-2}}(s\e-a))= exp_{q^{-2}}(1\otimes s\e-a)
exp_{q^{-2}}(s\e-a\otimes k_\a).
\label{sl2}
\ed
The relations (\ref{sl1}) and (\ref{sl2}) follow from the addition theorem for
$q$-exponents:
$$\exp_q(x+y)= \exp_qx\exp_qy\hsp{\mbox{if}}\;\;\; yx=qxy.$$ As a result
 we get the
 following commutation relations of $\Phi_\pm$ and $\Psi^\pm$ with
$\exp_{q^2}(t\ea )$ and $\exp_{q^{-2}}(s\e-a)$:
\bn
\exp_{q^2}(t\ea )\Phi_+=\Phi_+\exp_{q^2}(t\ea )
\label{8a}
\ed
\bn
\exp_{q^2}(t\ea )\Phi_-=(t\Phi_+k_\a^{-1}+\Phi_-)\exp_{q^2}(t\ea ),
\ed
\bn
\exp_{q^2}(t\ea )\Psi^+ =\Psi^+ \exp_{q^2}(qt\ea )-qt\Psi^-
\exp_{q^2}(q^{-1}t\ea ),
\ed
\bn
\exp_{q^2}(t\ea )\Psi^- =\Psi^-\exp_{q^2}(q^{-1}t\ea),
\ed
\bn
\Phi_+\exp_{q^{-2}}(s\e-a)=
\exp_{q^{-2}}(q^{-1}s\e-a)\Phi_+-\exp_{q^{-2}}(qs\e-a) q^{-1}s\Phi_-
\ed
\bn
\Phi_-\exp_{q^{-2}}(s\e-a)=\exp_{q^{-2}}(qs\e-a)\Phi_-,
\ed
\bn
\Psi^+\exp_{q^{-2}}(s\e-a)=\exp_{q^{-2}}(s\e-a)\Psi^+,
\ed
\bn
\Psi^-\exp_{q^{-2}}(s\e-a)=\exp_{q^{-2}}(s\e-a)(\Psi^-+sk_\a\Psi^+).
\label{8b}
\ed
The relations (\ref{8a})--(\ref{8b}) can be also deduced from
Proposition \ref{Kharchev}.

 We can use for instance $\Phi_\pm: V_{j-1/2} \rightarrow
V_j$ and $\Psi^\pm: V_{j'-1/2} \rightarrow V_{j'}$ and substitute
(\ref{8a})--(\ref{8b}) into basic equation (\ref{Q4}) for $$\G = \Phi
_+\otimes\Psi^+ +\Phi_-\otimes\Psi^- :  V_{j-1/2}\otimes
V_{j'-1/2}\rightarrow V_j\otimes V_{j'}.$$ Taking into account that $$
\Phi_+\mid j-1/2>=\mid j>,\hsp
\Phi_-\mid j-1/2>=\frac{q^{-2j+1}}{[2j]_q}\e-a\mid j>,$$
$$\Psi^+\mid j-1/2>=-\frac{q}{[2j]_q}\e-a\mid j>, \hsp
\Psi^-\mid j-1/2>=\mid j>,$$
 $$ <j\mid\Phi_+= <j-1/2\mid ,\hsp <j\mid\Phi_-=<j\mid\Psi^+=0,\hsp
 <j\mid\Psi^-= <j-1/2\mid$$
and
$$
\exp_{q^{-2}}(s\e-a)\e-a =\partial^{(q^{-2})}_s\exp_{q^{-2}}(s\e-a),
$$
 $$ \ea\exp_{q^2}(t\ea )= \partial^{(q^{2})}_t\exp_{q^{-2}}(t\ea)$$
where
\bn
\partial^{(q)}_x f(x)=\frac{f(qx)-f(x)}{(q-1)x}
\label{8c}
\ed
we get from (\ref{Q4}) the following Hirota identities:
$$
\left(
\frac{\partial_y^{(q^{-2})}}{[2j']_q}-\frac{q^{-2j}\partial_x^{(q^{-2})}}
{[2j]_q}+
(q^{2j'-2j-1}y-q^{2j-1}x)\frac{\partial_x^{(q^{-2})}}{[2j]_q}
\frac{\partial_y^{(q^{-2})}}{[2j']_q}\right)\tau_j(u,x)\tau_{j'}(v,y)=$$
\bn
=(v-q^{-2j}u)\tau_{j-1/2}(u,q^{-1}x)\tau_{j'-1/2}(q^{-1}v,y)
\label{lm}
\ed

We can use $q$-difference analog of Taylor formula
$$
f(x)=f(a)+(x-a)\frac{\partial^{(q)}}{\partial x}f(a)+\frac{(x-a)(x-qa)}{(2)_q!}
\frac{{\partial^{(q)}}^2}{\partial x^2}f(a)+\ldots +$$
\bn
\frac{(x-a)\cdots (x-q^{n-1}a)}{(n)_q!}
\frac{{\partial^{(q)}}^n}{\partial x^n}f(a)+\ldots +
\label{T}
\ed
and expand both sides of equation (\ref{lm}) in a series
\bn
\sum_{k,l\geq 0}P_{k,l}(x,u)(y-x')(y-q^{-2}x')\cdots (y-q^{-2(k-1)}x')
(v-u')(v-q^{-2}u')\cdots (v-q^{-2(l-1)}u')
\label{TTT1}
\ed
where $x'=xq^\a$, $u'=uq^\b$ for certain constants $\a$
 and $\b$ predicted
by the form of the equation and get as a result a hierarchy of bilinear
 $q$-difference equations
$$ P_{k,l}(x,u)=0$$
For instance, in a particular case $j=j'=\frac{1}{2}$ an equation (\ref{lm})
 looks like
\bn
\left( q\partial_y^{(q^{-2})}-\partial_x^{(q^{-2})}+(y-qx)
\partial_y^{(q^{-2})}\partial_x^{(q^{-2})}\right)
\tau_{1/2}(u,x)\tau_{1/2}(v,y)= qv-u
\label{1/2}
\ed
Here we fix the normalization $\tau_0=1$.

By definition, an equation (\ref{1/2}) has a solution
$$\tau_{1/2}(u,x)= a+bu+cx+dux$$
where $a= \tau_{<+\mid ,\mid +>}$, $b= \tau_{<+\mid ,\mid ->}$,
$c= \tau_{<-\mid ,\mid +>}$, $d= \tau_{<-\mid ,\mid ->}$. According to
 (\ref{RTT})  the matrix
$$T=\left(
\begin{array}{cc}
a&b\\ c&d
\end{array}
\right)$$
generates an algebra $Fun_q(SL_2)$ and satisfy the relations \cite{FRT}
$$
ab=q^{-1}ba,\;\; ac=q^{-1}ca,\;\; bd = q^{-1}db, \;\; cd = q^{-1}dc, \;\;
 bc=cb,$$
\bn
ad-q^{-1}bc = da- qbc=1
\label{F}
\ed
Note that the relation (\ref{1/2}) gives only half of the relations (\ref{F}).
 The rest of them one can deduce from the equation (\ref{Q4}) for vertex
 operators $\Phi_W^\pm: V_j \rightarrow V_{j-1/2}$,
$\Psi^W_\pm: V_{j'} \rightarrow V_{j'-1/2}$, $j= j' =1/2$.

An expansion (\ref{TTT1}) with $\a =1, \b =-1$ gives a hierarchy of
 equations with the first three listed below (for simplicity we write
in these equations $\tau$ instead of $\tau_{1/2}$).

\bn \tau(u,x)\cdot
\partial_x^{(q^{-2})}\tau(q^{-1}u,qx)=
\partial_x^{(q^{-2})}\tau(u,x)\cdot\tau(q^{-1}u,qx),
\label{1}
\ed
\bn
\tau(u,x)\cdot\partial_x^{(q^{-2})} \partial_u^{(q^{-2})} \tau(q^{-1}u,qx)-
 \partial_x^{(q^{-2})} \tau(u,x)\cdot\partial_u^{(q^{-2})} \tau(q^{-1}u,qx)=1
\label{2}
\ed
$$\tau(u,x)\cdot {\partial_{xx}^{(q^{-2})}}^2\tau(q^{-1}u,qx)-
\partial_x^{(q^{-2})} \tau(u,x)\cdot\partial_x
^{(q^{-2})}\tau(q^{-1}u,qx)+$$
\bn
+ q^2\partial_x^{(q^{-2})} \tau(u,x)\cdot
\partial_x^{(q^{-2})} \tau(q^{-1}u,q^{-1}x)=0
\label{3}
\ed
The first equation explaines the commutation rule of $\tau$ and its derivative,
 the second equation is a $q$-analog of Liouville equation and the third one
 should
restrict a solution to a simple bilinear form on $x$  .
 Differentiation of the first equation enables one to rewrite (\ref{2}) and
 (\ref{3}) in different forms; so we see that there is no assymmetry in $u,x$,
 $\tau$ and $\partial^2\tau$ in equations (\ref{2}) and (\ref{3}).

 In the classical limit $q=1$ $\tau_{u,x}$ is commutative and the
equation (\ref{1}) is trivial; the equation (\ref{2}) is
 $$\tau\partial_{ux}^2\tau-\partial_u\tau\partial_x\tau =1$$
which is Liouville equation in variable $\phi$, $\tau =e^{\phi}$:
$$ \phi_{xu}=e^\phi$$
and the equation (\ref{3}) restricts a solution $\tau$ to a linear form
 of $x$.

{\bf Remarks}. 1. Let us note once more that we treate $\tau$-function as
an element of algebra $Fun_q(G)$. In representation $g$ of $Fun_q(G)$
 $\tau_{u,x}$ turns to be an operator-valued function being a solution
 of corresponding equations (\ref{lm}), (\ref{1/2}),
(\ref{1})--(\ref{3}).\\
 2. A general solution $\tau(u,x)$ of the equations (\ref{lm}) is a
  generating function of matrix coeficients of matrix elements of  a
  representation of $U_q(sl_2)$ of spin $j$. It could be obtained
  directly from the universal $T$-matrix (see the next section).

 %%%%%%%%%%%%%%%%%%%%%%%%%%%%%%%%%%%%%%%%%%%%%%%
\setcounter{equation}{0}
\section{The universal $T$-matrix and group-like elements}
We can see that the ideology of intertwining operators allows one to write down
 a hierarchy of differential (or $q$-difference) equations provided there is
 known common bosonization of fundamental representations and of vertex
 operators $\Phi^W_i$, $\Phi_W^i$, $\Psi^W_i$ and $\Psi_W^i$, acting between
 them. The construction of such a bosonization is the main problem for the case
 of $U_q(\g)$.

However we can treate $\tau$-function as certain generating function of matrix
 elements and work directly with matrix elemens. We have then the basic
 relation (\ref{Q4}) defining the connection between matrix elements of
 different representations. There is also a canonical way for producing matrix
 elements provided the structure of representation is known. Let us restrict
 ourselves to a case of $\Aa =U_q(\g)$ where $\g$ is simple finite-dimensional
Lie  algebra.

Let $H=\Aa\otimes \A* = U_q(\g)\otimes Fun_q(G)$ and $T$ is a canonical tensor
 of a Hopf pairing of $\Aa$ and $\A*$:
$$ T=\sum e_i\otimes e^i,\hsp \{ e_i,e^j\} =\delta_{i,j}.$$
 $T$ is usually called as universal $T$-matrix \cite {FRT}, \cite{FG}.
 The properties of Hopf pairing in terms of $T$ mean that
\bn
(\D\otimes id)T=T^{13}T^{23},
\label{T1}
\ed
\bn
(id\otimes\D)T=T^{12}T^{13},
\label{T2}
\ed
\bn
(S\otimes id)T=T^{-1},
\label{T3}
\ed
\bn
(id\otimes S)T=T^{-1}.
\label{T4}
\ed

Let $\rho :\Aa\rightarrow {\mbox{End}}\, V$ be a representation of
$\Aa =U_q(\g)$ in a space $V$ with basis $\{ v_i\}$. Then  matrix elements
 $ <v^j\mid (\rho\otimes id)T\mid v_i>$
 are, by definition of $T$ , the elements $\tau_{v^j,v_i}$ of an algebra
 $\A* =Fun_q(G)$:
 $$\tau_{v^j,v_i}= <v^j\mid (\rho\otimes id)T\mid v_i>.$$
 Furtermore, for a representation $g$: $\A* =Fun_q(G) \rightarrow
 {\mbox{End}}\, U$ an operator $\tau_{v^j,v_i}^g$ can be also expressed
 in terms of $T$:
 $$\tau_{v^j,v_i}^g= <v^j\mid (\rho\otimes g)T\mid v_i>.$$
 The universal $T$-matrix can be written in factorized form
 (see \cite{FG} for $U_q(gl_n)$ case):
\bn
T=\prod_{\gamma\in \D_+}^{\rightarrow}\exp_{q^{-(\gamma,\gamma)}}
a(\gamma)e_\gamma\otimes s_\gamma\cdot
q^{d_{i,j}h_i\otimes \hat{h}_j}\cdot
\prod_{\gamma\in \D_+}^{\leftarrow}\exp_{q^{(\gamma,\gamma)}}
-a(\gamma)e_{-\gamma}\otimes s_{-\gamma}
\label{TT}
\ed
 Here
$e_{\pm\gamma}\in U_q(\bb_\pm)$ are Cartan--Weyl generators of $U_q(\g)$
 constructed for a fixed reduced decomposition of the longest element
 $\omega_0$ of $q$-Weyl group or, equivalently, for a fixed normal
ordering of a system of positive roots $\D_+$; $a(\gamma)$ are the
 constant coeffitients normalizing the relation
$$[e_{\gamma},e_{-\gamma}]=\frac{k_\gamma-k_\gamma^{-1}}{a(\gamma)}$$
 and $d_{i,j}$ is an inverse matrix to symmetrized Cartan matrix of $\g$.
 $U_q(\bb_\pm)$ are Borel subalgebras of $U_q(\g)$.

 The elements $s_\gamma,  \gamma\in\D_+$ satisfy the same commutation
  relations as $e_{-\gamma}\in U_q(\bb_-)$,
 the elements $s_{-\gamma},  \gamma\in\D_+$ satisfy the same commutation
  relations as $e_{\gamma}\in U_q(\bb_+)$,
   $q^{\hat{h}_i} s_{\gamma}= q^{(\a_i,\gamma)} s_{\gamma}
    q^{\hat{h}_i}$,
   $q^{\hat{h}_i} s_{-\gamma}= q^{(\a_i,\gamma)} s_{-\gamma}
    q^{\hat{h}_i}$ and $[s_\gamma ,s_{-\delta}] =0$ for any
  $\gamma  ,\delta \in  \D_+$.

 An expression (\ref{TT}) can be deduced from analogous expression for
 the universal $R$-matrix for $U_q(\g)$ \cite{KR}, \cite{KT}, \cite{LS}.

 %%%%%%%%%%%%%%%%%%%%%%%%%%%%%%%%%%%%%%%%%%%
To show this we consider the double ${\cal D}={\cal D}(U_q({\g}))$ of
 $U_q({\g})$.  The universal  $R$-matrix for ${\cal D}$ belongs to $(
U_q({\g})\otimes  1) \otimes (1 \otimes Fun_q^0(G))$ where ${}^0$
denotes an opposite comultiplication and enjoys the properties:  \bn
\Delta'(x)R=R\D(x)
\label{R1}
\ed
\bn
(\D\otimes id)R=R^{13}R^{23},
\label{R2}
\ed
\bn
(id\otimes\D)R=R^{13}R^{12},
\label{R3}
\ed
\bn
(S\otimes id)R=R^{-1},
\label{R4}
\ed
\bn
(id\otimes S')R=R^{-1}.
\label{R5}
\ed
We are not interested here in the property (\ref{R1}) and so have no need in
 commutation relations between $U_q(\g)$ and $Fun_q^0(G)$ inside
${\cal D}(U_q({\g}))$. In this setting we can identify $U_q(\g)$ with
$( U_q{\g}\otimes  1)$ and $Fun_q(G)$ with $(1 \otimes Fun_q^0(G))$ returning
 to original comultiplication. Under these identifications the universal
$R$-matrix for ${\cal D}(U_q({\g}))$ coincides with universal $T$-matrix $T$,
 the properties (\ref{R2})--(\ref{R5}) correspond to (\ref{T1})--(\ref{T4}).

 We reconstruct   the universal $R$-matrix for
${\cal D}(U_q({\g}))$ from the universal $R$-matrix
 for $U_q(\g)$  using the isomorphisms of algebras
\bn
Fun_q(G)\simeq U_q(\bb_+)\otimes_{U(\h)}^{}U_q(\bb_-)
\label{iso}
\ed
where  $U_q(\bb_\pm)$ are Borel subalgebras of $U_q(\g)$ and ${U(\h)}$
is Cartan
 subalgebra.

Let $R_U$ be the universal $R$-matrix for $U_q(\g)$. In the same notations
 it
has a factorized form \cite{KR}, \cite{LS}, \cite{KT}:  $$ R_U=R\cdot
K,$$ where \bn R=\prod_{\gamma\in
\D_+}^{\rightarrow}\exp_{q^{-(\gamma,\gamma)}} a(\gamma)e_\gamma\otimes
e_{-\gamma}, \hsp K=q^{d_{i,j}h_i\otimes h_j}; \label{URM} \ed
 The tensor $R_U$
 defines a pairing of $U_q(\bb_+)$ with  $U_q(\bb_-)$.

We can
 rewrite its component $R$ as \bn R=\prod_{\gamma\in
\D_+}^{\rightarrow}\exp_{q^{-(\gamma,\gamma)}}a(\gamma) e_\gamma\otimes
s_\gamma
\label{R}
\ed
where $s_\gamma =  e_{-\gamma}$ is dual to $a(\gamma)e_\gamma$ in $U_q(\bb_-)$.
 Analogously, $(R_U^{21})^{-1}$ is again the universal $R$-matrix for
 $U_q(\g)$. It looks like
$$(R_U^{21})^{-1}= K^{-1}\overline{R}$$
where
$$\overline{R}=\prod_{\gamma\in \D_+}^{\leftarrow}\exp_{q^{(\gamma,\gamma)}}
-a(\gamma)e_{-\gamma}\otimes s_{-\gamma}$$
with $s_{-\gamma} =  e_{\gamma}$ being dual to $-a(\gamma)e_{-\gamma}$ in
$U_q(\bb_+)$.

The decomposition (\ref{iso}) and the fact, that  $U_q(\bb_\pm)$ are Hopf
 subalgebras of $U_q(\g)$ mean that we can compose the pairing of $U_q(\g)$
 with $Fun_q(G)$ from $R_U$ and $(R_U^{21})^{-1}$ after proper factorization
 over $U(\h)$, in other words,
$$T=RK\overline{R},$$ and we get the expression(\ref{TT}).

 Since $s_{\pm\gamma}$ and $\hat{h}_j$ are supposed to commute with $U_q(\g)$
 we can treate $T$ as a product of $q$-exponents of $e_{\pm\gamma}$ with
 parameters $s_{\pm\gamma}$ and $\hat{h}_j$:
\bn
T(s_{\pm\gamma},\hat{h}_j)=
\prod_{\gamma\in \D_+}^{\rightarrow}\exp_{q^{-(\gamma,\gamma)}}
e_\gamma s_\gamma\cdot
q^{d_{i,j}h_i \hat{h}_j}\cdot
\prod_{\gamma\in \D_+}^{\leftarrow}\exp_{q^{(\gamma,\gamma)}}
e_{-\gamma} s_{-\gamma}
\label{TT1}
\ed
 The property (\ref{T1}) applied to a representation $\rho\otimes \rho'$ of
 $U_q(\g)\otimes U_q(\g)$ means that
$$ \rho\otimes \rho'\D T(s_{\pm\gamma},\hat{h}_j)=
\rho (T(s_{\pm\gamma},\hat{h}_j)\otimes \rho'T(s_{\pm\gamma},\hat{h}_j).$$
It could be treated as ''group-like'' property of $T$. The property (\ref{T2})
 means that
$$ T(s_{\pm\gamma}^1,\hat{h}_j^1)T(s_{\pm\gamma}^2,\hat{h}_j^2)=
T(s_{\pm\gamma}^3,\hat{h}_j^3)$$
for mutually commuting parameters $\{ s_{\pm\gamma}^1,\hat{h}_j^1\}$,
$\{ s_{\pm\gamma}^2,\hat{h}_j^2\}$; where parameters $\{
 s_{\pm\gamma}^3,\hat{h}_j^3\}$ satisfy the same commutation relations as
$\{ s_{\pm\gamma}^1,\hat{h}_j^1\}$ and $\{ s_{\pm\gamma}^2,\hat{h}_j^2\}$ by
 themselves. Analogously, the last property (\ref{T4}) means that
$T^{-1}(s_{\pm\gamma},\hat{h}_j)$ is group-like element for an opposite
 comultiplication $\D'$.
%%= T({s'}_{\pm\gamma},{\hat{h}'}_j)$ for some
%% other $\{ {s'}_{\pm\gamma},{\hat{h}'}_j\}$.

This reformulation gives a hint to treate $T$ in (\ref{TT1}) as a $q$-analogue
 of a parametrization of Chevalley group \cite{S}. It is possibble to try to
 follow  the games which usually people play with Chevalley groups:

 One can define this ''group'' as a group generated by $q$-exponents of
 simple root vectors. Such a parametrization was presented in \cite{MV} in a
 pretty nice form. It is also possible to write down the defining relations
 for the generators of this ''group''. The relations turn to be $q$-analogs
 of classical Chevalley relations. First they appear in \cite{KT}, eq.
(5.10)-(5.12)
and then in \cite{MV} in the context of universal $T$-matrix.
\vspace{.5cm}

{\large{\bf Acknowlegements. }} The authors are grateful to A. Gerasimov,
 A. Mironov, A. Morozov and A. Zabrodin for numerous helpful
 discussions. D.L. thanks Uppsala University for kind hospitality.
S.Khor. was partially supported by ISF grant MBI000,
S.Khar. was partially supported by grant 93-02-14365 of the Russian
Foundation of Fundamental Research and by ISF grant MGK000, and D.L.
was partially supported by ISF grant MIF000, AMS FSU foundation and
by grant 93-02-14365 of the Russian Foundation of Fundamental Research.

%%%%%%%%Subject: hirota.bib


\begin{thebibliography}{99}

\bibitem[DJ1]{D}
 Date, E., Jimbo, M., Kashiwara, M. and Miwa, T. Transformation
 groups for soliton equations, in ''Nonlinear integrable systems -
		  classical theory and quantum theory'' {\it World
		  Scientific} (1983).
\bibitem[DJ2]{XXZ}
Date, E., Jimbo, M., Kashiwara, M., Miwa, T. and Nakayashiki, A.
 Diagonalization of the $XXZ$ Hamiltonian by Vertex Operators.
 {\it Comm. Math. Phys., 151} (1993), 89-153.

\bibitem[Dr]{Dr}
Drinfeld, V.G. Quantum groups.
{\it Proc. ICM-86 (Berkely USA) vol.1},
		  798-820. Amer. Math. Soc. (1987).

\bibitem[F]{F}
Faddeev, L.,D. Quantum completely integrable models in field theory.
 {\it Sov. Sci. Rev. Math. Phys. C1} (1980), 107-155.

\bibitem[FRT]{FRT}
Faddeev, L.D., Reshetikhin, N.Yu. and Takhtajan, L.A. Quantization of
		  Lie groups and Lie algebras. {\it Algeebra and
		  Analisys 1} (1989), 178-201.

\bibitem[FR]{FR}
Frenkel, I.B., and Reshetikhin, N.Yu.
Quantum Affine  Algebras and Holonomic Difference equations.
{\it Commun. Math. Phys., 146} (1992), 1-60.

\bibitem[FG]{FG}
 Fronsdal, Galindo. The Universal $T$-matrix, {\it preprint UCLA/93/TEP/2}
 (1993).

\bibitem[GKLMM]{GKLMM}
 Gerasimov, A., Khoroshkin, S., Lebedev, D., Mironov, A., Morozov, A.
 Generalized Hirota Equations and Representation Theory.
 {\it ITEP M-2/94, hep-th/9405011} (1994).

\bibitem[GS]{GS}
 Gervais, J-L and Schnittger, J. The many faces of quantum Liouville
 exponentials. {\it LPTENS-93/30, hep-th/9308134}.

\bibitem[JM]{JM}
 Jimbo, M. and Miwa, T. Algebraic Analisys of Solvable Lattice Models.

\bibitem[Kac]{Kac}
 Kac, V., Infinite-dimensional Lie algebras. {\it Cambridge
 Univ. Press} (1985).

\bibitem[KR]{KR}
Kirillov, A.N and Reshetikhin, N.Yu.
 q-Weyl group and a
multiplicative formula for universal $R$-matrices.
{\it Comm. Math. Phys. 134} (1990), 421-431.

\bibitem[KT]{KT}
Khoroshkin, S.M., and Tolstoy, V.N.
Universal R-matrix for quantized (su\-per)\-algebras.
{\it Commun. Math. Phys. 141} (1991), 599-617.

\bibitem[KW]{KW}
Kac, V.G., Wakimoto, M. Exceptional hierarchies of Solitoon equations.
 {\it Proc. Symposia in Pure Math., 49} (1989), 191-237.

\bibitem[LS]{LS}
Levendorskii, S.Z., and Soibelman, Ya.S. Some application of
quantum Weyl groups. The multiplicative formula for  universal
$R$-matrix for simple Lie algebras. {\it J. Geom. and Phys 7:4} (1990),
1-14.

\bibitem[MV]{MV}
Morozov, A., and Vinet, L.
Free-Field Representation of Group Element for Simple Quantum Group.
 {\it ITEP-M3/94 and CRM-2202 preprints, hep-th/9409093} (1994).


\bibitem[S]{S}
Steinberg, S., Lectures on Chevalley Groups. {\it Yale Univ. Press}
 (1967).


\bibitem[UT]{UT}
Ueno, K., Takasaki, K. Toda lattice hierarchy. {\it Adv. Studies in
 Pure Math.} (1984).

\bibitem[ZZ]{ZZ}
Zamolodchikov, A.B., Zamolodchikov, Al.B., Factorized $S$-matrix in two
 dimensions as the exact solutions of certain relativistic quantum field
 theory models. {\it Ann. Phys. 120} (1979), 253-291.
\end{thebibliography}
\end{document}